\documentclass[reqno]{alea2}
\usepackage{natbib}
\usepackage{fancyhdr}
\usepackage{graphicx}

\usepackage{amssymb}
\usepackage{amsmath}

\renewcommand{\cite}{\citet}

\makeatletter \@addtoreset{equation}{section} \makeatother

\renewcommand\theequation{\thesection.\arabic{equation}}
\renewcommand\thefigure{\thesection.\@arabic\c@figure}
\renewcommand\thetable{\thesection.\@arabic\c@table}

\theoremstyle{plain}
\newtheorem{theorem}{Theorem}[section]

\theoremstyle{definition}

\theoremstyle{remark}

\newtheorem{example}[theorem]{Example}

\def\cC{\mathcal{C}}

\newtheorem{theo}{Theorem}[section]
\newtheorem{prop}[theo]{Proposition}
\newtheorem{corol}[theo]{Corollary}
\newtheorem{lemme}[theo]{Lemma}
\newtheorem{Rem}[theo]{Remark}

\newcommand{\CQFD}{\hfill $\square$}

\newcommand{\ind}{\mathbf{1}}

\renewcommand{\MakeUppercase}[1]{#1}

\def\real{\mathbb{R}}

\def\nit{\Bbb{N}}
\def\pit{\Bbb{P}}

\def\ee{\Bbb{E}}

\def\pit{\Bbb{P}}
\def\E{\mathop{\hbox{\rm I\kern-0.20em E}}\nolimits}
\def\Var{\mathop{\hbox{\rm Var}}\nolimits}

\def\og{\leavevmode\raise.3ex
     \hbox{$\scriptscriptstyle\langle\!\langle$~}}
\def\fg{\leavevmode\raise.3ex
     \hbox{~$\!\scriptscriptstyle\,\rangle\!\rangle$}~}

\newcommand{\cD}{{\mathcal D}}

\newcommand{\cM}{{\mathcal M}}

\newfont{\msbm}{msbm10 scaled\magstep1}
\newfont{\msbms}{msbm7 scaled\magstep1} 

\newcommand{\bbP}{\mbox{$\mbox{\msbm P}$}}

\begin{document}


\author{ Jean-Christophe Breton}
\address{Laboratoire MIA, Universit\'e de La Rochelle, 17042 La Rochelle Cedex, France }
\email{jcbreton@univ-lr.fr}
\author{Cl\'ement Dombry}
\address{Laboratoire LMA, Universit\'e de Poitiers, T\'el\'eport 2, BP 30179, F-86962 Futuroscope-Chasseneuil cedex, France}
\email{clement.dombry@math.univ-poitiers.fr}

\title[Weighted random ball model]{Functional macroscopic behavior of weighted random ball model} 

\begin{abstract}
We consider a generalization of the weighted random ball model defined by 
$$
M({\bf y})=\int_{\real^d\times\real^+\times\real}mh\left(\frac{{\bf y}-{\bf x}}{r}\right)N(d{\bf x},dr,dm)
$$
where $N$ is a random Poisson measure on $\real^d\times\real^+\times\real$ with a product heavy tailed intensity measure 
and $h:\real^d\to \real$ is a fading function. 
This functional can serve as a basic model for transmission with fading effect. 
In \cite{BD}, the authors proved the convergence of the finite-dimensional distributions of related generalized random fields under various scalings and in the particular case when $h$ is the indicator function of the unit ball in $\real^d$. 
In this paper, tightness and functional convergence are investigated. 
Using suitable moment estimates, we prove functional convergences for some parametric classes of configurations under the so-called large ball scaling and intermediate ball scaling. 
Convergence in the space of distributions is also discussed.
\end{abstract}

\maketitle

\section{Introduction}
\label{sec:intro}
 
We consider weighted random balls in $\real^d$ generated by a Poisson random measure $N_\lambda$ on $\real^d\times\real^+\times\real$  with intensity 
$$
n_\lambda(d{\bf x},dr,dm)=\lambda d{\bf x} F(dr) G(dm)
$$
where $\lambda\in\real^+$ and $F$, $G$ are probability measures on $\real^+$ and $\real$  respectively. 
For each $3$-uplet $({\bf x},r,m)$, ${\bf x}$ represents the center of the Euclidean ball $B({\bf x},r)$ and $r$ its radius, 
$m$ stands for the weight of the ball. 
The parameter $\lambda$ is interpreted as the intensity of the balls in $\real^d$.
Such models are used for instance to represent a spatial communication network, see \cite{Kaj-net}, \cite{YP}. 
In this case, ${\bf x}$ represents a station transmitting a signal, $r$ the range of emission and $m$ the intensity of the signal, 
see \cite{BD} and reference therein in particular \cite{KLNS}, \cite{KT07}. 
Following this interpretation, the signal $m$ transmitted by ${\bf x}$ is received in some ${\bf y}\in\real^d$ if and only if ${\bf y}\in B({\bf x},r)$,  
and the quantity of signal received from the stations at ${\bf y}$ is given by
\begin{equation}
\label{eq:B}
\int_{\real^d\times\real^+\times\real} m\ind_{B({\bf x},r)}({\bf y}) N_\lambda(d{\bf x},dr,dm)=\int_{\real^d\times\real^+\times\real} m\ind_{B({\bf 0},1)}(({\bf y}-{\bf x})/r) N_\lambda(d{\bf x},dr,dm).
\end{equation}
From a modeling point of view, it is natural to consider that the signal transmitted by ${\bf x}$ and received in ${\bf y}$ fades when ${\bf y}$ gets away from the station ${\bf x}$. 
In order to take into account this phenomenon, we introduce a fading function $h$ replacing $\ind_{B({\bf 0},1)}$ in \eqref{eq:B}, 
more precisely the faded signal received at ${\bf y}$ from ${\bf x}$ is $mh(({\bf y}-{\bf x})/r)$ 
and, assuming moreover that no interference occurs between the stations, 
the quantity of signal received at ${\bf y}$ is now given  by
\begin{equation}
\label{eq:h}
M({\bf y})=\int_{\real^d\times\real^+\times\real} mh(({\bf y}-{\bf x})/r) N_\lambda(d{\bf x},dr,dm).
\end{equation}
From a physical point of view, it is natural to assume that $h({\bf y})$ is a radially non-increasing function with $h({\bf 0})=1$, $0\leq h({\bf y})\leq 1$ and $\lim_{\|{\bf y}\|\to+\infty} h({\bf y})=0$. 
The function $h$ is said to be radially non increasing if for all ${\bf y}\in\real^d$, the function $r\mapsto h(r{\bf y})$ is non-increasing on $[0,+\infty)$. 
However, on a mathematical point of view, more general assumptions will be enough, 
see assumption \eqref{eq:halpha} below.

\medskip
In the sequel, we are more generally interested in the contribution 
$$M(\mu)=\int_{\real^d} M({\bf y})\mu(d{\bf y})$$ 
of the model in a configuration of points ${\bf y}$ represented by a measure $\mu$. 
For instance, the configuration reduced to the point ${\bf y}$ is represented by $\mu=\delta_{{\bf y}}$ 
and in this case $M(\delta_{{\bf y}})=M({\bf y})$.
It is natural to consider finite positive measures $\mu$ on $\real^d$ but our study supports signed measures $\mu$ with finite total variation. 
In the sequel, we shall note $\mathcal{M}$ the set of such measures and we recall that, 
equipped with the total variation norm $\|\mu\|_{\mathcal M}=|\mu|(\mathbb R^d)$, 
$\mathcal{M}$ is a Banach space. Actually in order to make our study easier and in contrast with \cite{BD}, we shall consider signed measures $\mu$ with density, {\it i.e.} $\mu(d{\bf y})=\phi({\bf y})d{\bf y}$ for some $\phi\in L^1(\real^d)$. Setting $\tau_{{\bf x},r}h({\bf y})=h(({\bf y}-{\bf x})/r)$ and $\mu[f]=\int_{\real^d} f({\bf y})\mu(d{\bf y})$ for $f\in L^1(\real^d,\mu)$, 
the Fubini theorem allows to rewrite  
\begin{equation}
\label{eq:M}
M(\mu)=\int_{\real^d\times\real^+\times\real} m\mu[\tau_{{\bf x},r}h ] N_\lambda(d{\bf x},dr,dm).
\end{equation}
Note that the stochastic integral in \eqref{eq:M} is well defined 
and the change in the order of integrals is authorized when
\begin{eqnarray}
\label{eq:Fub}
\nonumber
&&\int_{\real^d\times\real^+\times\real} |m\mu[\tau_{{\bf x},r}h ]| \ n_\lambda(d{\bf x},dr,dm)\\
&\leq &\lambda |\mu|(\real^d) \left(\int_{\real^d}|h({\bf x})|d{\bf x}\right)\left(\int_{\real^+} r^d F(dr)\right)\left(\int_{\real} |m| G(dm)\right) <+\infty. 
\end{eqnarray}
We will always suppose that all three integrals in \eqref{eq:Fub} above are finite (precise assumptions on $F$, $G$ and $h$ are given in the set of conditions $(\bf A)$ below).
Furthermore, in this case, the expected value of $M(\mu)$ is given by
$$
\ee[M(\mu)]= \lambda \mu(\real^d) \left(\int_{\real^d}h({\bf y})d{\bf y}\right)  \left(\int_{\real^+} r^d F(dr)\right)\left(\int_{\real} m G(dm)\right) .
$$
 \ \\
In order to investigate the macroscopic behavior of the generalized random field $(M(\mu))_{\mu\in{\mathcal M}}$, 
we apply a scaling $x\mapsto \rho x$, with $\rho<1$. 
The scaling contracts the space $\real^d$ and is interpreted as zoom-out in the model.
Note that when $\rho>1$, the scaling becomes a dilatation of ${\mathbb R}^d$ and is interpreted as zoom-in. 
In contrast to  \cite{BEK}, \cite{BD}, we focus in this note only on zoom-out (see below for further comments on the relation between this note and our previous paper \cite{BD}). 
In order to derive non-trivial asymptotics, the intensity $\lambda$ of the Poisson measure is adapted to the scaling procedure by allowing $\lambda:=\lambda(\rho)$ to depend on the zooming factor $\rho$, see \cite{BD}. 
Note that the natural intensity $\lambda(\rho)$ corresponding to the scaling $x\mapsto \rho x$ is $\lambda(\rho)=\rho^{-d}\lambda$.
After the scaling, the generalized field becomes 
\begin{equation}
\label{eq:Mrho}
M_{\rho}(\mu)=\int_{\real^d\times\real^+\times\real} m \mu[\tau_{{\bf x},r}h ] N_{\rho,\lambda(\rho)}(d{\bf x},dr,dm)
\end{equation}
where $N_{\rho,\lambda(\rho)}$ is the Poisson random measure with intensity
$$
n_{\rho,\lambda(\rho)}(d{\bf x},dr,dm)=\lambda(\rho)d{\bf x}F_\rho(dr)G(dm)
$$  
and $F_\rho$ is the image measure of $F$ under $r\mapsto\rho r$. 
We are finally led to investigate, for a proper normalization $n(\rho)$, the fluctuations of the rescaled and centered random field 
$$
\widetilde M_\rho(\mu)=n(\rho)^{-1}(M_\rho(\mu)-\ee[M_\rho(\mu)]).
$$
\ \\
In order to derive non-trivial asymptotics for the model \eqref{eq:Mrho}, the  distributions $F$ and $G$ driving the behavior of the radius $r$ and of the weights $m$, 
and the shape function $h$ must satisfy some conditions, denoted by conditions $(\bf A)$, 
that we state now precisely:
\begin{itemize}
 \item The probability $G$ is assumed to belong to the normal domain of attraction of the $\alpha$-stable distribution $S_\alpha(\sigma,b,\tau)$ with $\alpha\in (1,2]$, {\it i.e.} if $X_1, \dots, X_n$ are {\it i.i.d.} according to $G$, $n^{-1/\alpha}(X_1+\dots+X_n)\Rightarrow  S_\alpha(\sigma,b,\tau)$. 
According to \cite[XVII.5]{Feller}, this is equivalent to the following estimate on the characteristic function $\varphi_G$ of~$G$: 
\begin{equation}
\label{eq:local}
\tag{${\bf A}_1$}\varphi_G(\theta)=1+i\theta\tau-\sigma^\alpha|\theta|^\alpha(1+i b \varepsilon(\theta)\tan(\pi\alpha/2))+o(|\theta|^\alpha) \quad \mbox{as}\ \theta\to 0.
\end{equation}
In the case $\alpha\in (1,2)$, a typical choice for $G$ is an heavy-tailed distribution while for $\alpha=2$, $G$ may be any distribution with finite variance. Observe that since $\alpha>1$, $G$ has a finite moment of order $1$.
\item The probability $F$ is assumed to have a regularly varying tail satisfying
\begin{equation}
\label{eq:F2}\tag{${\bf A}_2$}
\bar F(r):=\int_r^{+\infty} F(du)\sim_{r\to +\infty} C_\beta r^{-\beta}\quad \mbox{for some } d<\beta<\alpha d. 
\end{equation}
Here and in the sequel, $f(r)\sim_{r\to+\infty} g(r)$ indicates that $\lim_{r\to+\infty} f(r)/g(r)=1$.
Observe that, under \eqref{eq:F2}, the expectation of the volume of the random balls $\int_{\real^+} r^dF(dr)$ is finite 
(see Lemma \ref{lemme:F_truncated} below) and the bound \eqref{eq:Fub} indeed holds true.   
\item The shape function $h$ is assumed to be continuous almost everywhere and such that 
\begin{equation}
\label{eq:halpha}\tag{${\bf A}_3$}
 h^\ast({\bf x}):=\sup\{|h(r{\bf x})|\ : \ r\geq 1\}
 \in L^1(\real^d)\cap L^\alpha(\real^d).
\end{equation} 
Note that this implies that $h\in L^1(\real^d)\cap L^\alpha(\real^d)$ and that  if $h$ is radially non-increasing, then $h^\ast=h$. 
Indeed, $h^\ast$ is the smallest radially non-increasing function dominating $h$.
\end{itemize}

\medskip
The convergences of the finite-dimensional distributions ({\it fdd}) of $\widetilde M_\rho$ were (essentially) already derived in \cite{BD} 
under three different regimes depending on the behavior of $\lambda(\rho)\rho^\beta$ (the so-called large, small and intermediate ball regimes). 
In this note, we actually focus on the corresponding functional convergence for the generalized random fields $\left(\widetilde M_\rho(\mu)\right)_{\mu\in\cM}$. 
Since there is no natural functional space in which the random function $\mu\in{\mathcal M}\mapsto\widetilde M_\rho(\mu)$ belongs (at least heuristically), we choose to consider a special parametric sub-family $(\mu_{\bf t})_{{\bf t}\in\real^p}$ of ${\mathcal M}$ 
and to investigate the tightness of the random fields $\left(\widetilde M_\rho(\mu_{\bf t})\right)_{{\bf t}\in\real^p}$. 
Our main contribution is to prove tightness in the space of continuous functions $\cC(\real^p)$ of such random fields under suitable conditions on the family $(\mu_{\bf t})_{{\bf t}\in\real^p}$. 
As a by-product, we also obtain weak convergence of the random fields $\widetilde M_\rho$ in the space of distributions.

\medskip
Let us conclude this introduction with some comments on the relations between this note, our previous paper \cite{BD} and other related papers \cite{MRRS}, \cite{KT07}.
In \cite{BD}, {\it fdd} convergences are obtained for $(\widetilde M_\rho(\mu))_\mu$ when $\mu$ belongs to some special subspace ${\mathcal M}_{\alpha,\beta}$ on which 
\begin{equation}
\label{eq:Malphabeta}
\int_{\real^d}|\mu(B({\bf x},r))|^\alpha d{\bf x} \leq C(r^p\wedge r^q)\quad \mbox{ for some } p<\beta<q
\end{equation}
(roughly speaking the condition requires a control of the measures $\mu(B({\bf x},r))$ of both large and small balls, uniform in the centers of the balls). 
Moreover, \cite{BD} deals simultaneously with the macroscopic behavior ({\it i.e.} $\rho\to 0$ and $F$ has a power law behavior in $+\infty$ of order $\beta>d$) 
and microscopic behavior ({\it i.e.} $\rho\to+\infty$ and $F$ has a power law behavior in $0$ of order $\beta<d$). 
In comparison, in this note we deal only with the macroscopic behavior ({\it i.e.} $\rho\to 0$ and $\beta>d$) and for special measures $\mu(d{\bf y})=\phi({\bf y})d{\bf y}\in{\mathcal M}$ with density $\phi\in L^1(\real^d)\cap L^\alpha(\real^d)$. 
But, $F$ is not any more assumed to have a density and we only assume the tail condition \eqref{eq:F2}. 
Moreover a shape function $h$ is added in the model \eqref{eq:h} to take into account the fading of the signal in the communication network. 
But our main contribution in this setting is to derive tightness to obtain functional counterparts of the convergences in \cite{BD}.
In the particular case of the dimension $d=1$, the model is related to the infinite Poisson model in \cite{MRRS} and to the continuous flow reward in \cite{KT07}. 
In both papers, the authors deal with the asymptotic behaviour of the random process corresponding in our setting to  $(\widetilde M_\rho(1_{[0,t]}( y)d y))_{t\geq 0}$ and the issue of tightness is adressed in both papers.

\medskip
The rest of the note is organized as follows: 
the main results are stated in Section~\ref{sec:result} and proved in Section \ref{sec:proof}.  
Technical results are postponed in an Appendix. 


\section{Main results}
\label{sec:result}
First, we recall the finite-dimensional convergence for the generalized random field $\widetilde M_\rho=n(\rho)^{-1}(M_\rho(\mu)-\ee[M_\rho(\mu)])$. 
Actually, we state a slight modification of the main results in \cite{BD} replacing $\ind_{B({\bf 0},1)}$ therein by a shape function $h$ satisfying \eqref{eq:halpha}.
We shall abusively write $\mu\in L^1(\real^d)\cap L^\alpha(\real^d)$, instead of $\mu(d{\bf y})=\phi({\bf y})d{\bf y}$ with $\phi\in L^1(\real^d)\cap L^\alpha(\real^d)$.
The symbol $\stackrel{{\it fdd}}{\Longrightarrow}$ stands for the {\it fdd} convergences and the symbol $\stackrel{{\mathcal X}}{\Longrightarrow}$ is also used throughout 
to indicate a functional convergence in the functional space ${\mathcal X}$ (for instance, in Theorem \ref{theo:CV_functional1}, ${\mathcal X}={\mathcal C}(\real^d)$, the space of continuous function of $\real^d$). 

\begin{prop}
\label{prop:fdd} 
Suppose conditions $(\bf A)$ hold. 
\begin{enumerate}
\item (Large ball regime) 
If $\lambda(\rho)\rho^\beta\to+\infty$, then, setting $n(\rho)=(\lambda(\rho)\rho^\beta)^{1/\alpha}$, we have as $\rho\to 0$: 
\begin{equation}
\nonumber
\widetilde M_\rho(\mu) \stackrel{{\it fdd}}{\Longrightarrow} Z_\alpha(\mu),\quad \mu\in L^1(\real^d)\cap L^\alpha(\real^d)
\end{equation}
where  $Z_\alpha$ is the stable field 
$$
Z_\alpha(\mu)=\int_{\real^d\times\real^+} \mu[\tau_{{\bf x},r}h]M_\alpha(dr,d{\bf x})
$$
with respect to the $\alpha$-stable measure $M_\alpha$ with control measure
$\sigma^\alpha C_\beta r^{-1-\beta} drd{\bf x}$ and constant skewness function $b$, where $\sigma$ and $b$ are related to $G$ by \eqref{eq:local}. 
\item (Intermediate ball regime)
If $\lambda(\rho)\rho^\beta \to a$ for some $a\in (0,+\infty)$, then, setting $n(\rho)=1$, we have as $\rho\to 0$: 
$$
\widetilde M_\rho(\mu) \stackrel{{\it fdd}}{\Longrightarrow} J_a(\mu),\quad \mu\in L^1(\real^d)\cap L^\alpha(\real^d)
$$
where $J_a$ is the compensated Poisson integral
\begin{equation*}
\label{eq:J}
J_a(\mu)=\int_{\real^d\times\real\times\real^+} m\mu[\tau_{{\bf x},r}h] \widetilde N_{\beta,a}(d{\bf x},dr,dm)
\end{equation*}
with respect to the compensated Poisson random measure $\widetilde N_{\beta,a} $ with intensity  $aC_\beta r^{-\beta-1}d{\bf x}drG(dm)$.
%
\end{enumerate}
\end{prop}
In the sequel, the finite-dimensional results are strengthened into functional convergence for a parametric sub-family of measures $\mu_{\bf t}(d{\bf y})=\phi_{\bf t}({\bf y}) d{\bf y}, {\bf t}\in\real^p$, in $L^1(\real^d)\cap L^\alpha(\real^d)$, 
so that the generalized random field $\widetilde M_\rho$ induces a $p$-dimensional random field $(\widetilde M_\rho(\mu_{\bf t}))_{{\bf t}\in\real^p}$. 
To that aim, we investigate the tightness in ${\mathcal C}(\real^p)$ of this induced $p$-dimensional random field. 
The proof relies on a Censov criterion and on moment estimates for increments presented in Section \ref{sec:moment}. 
Actually in our setting, the relevant increments are generalized increments on blocks defined as follows. 
Let ${\bf s}, {\bf t}\in\real^p$ be such that $s_i\leq t_i, 1\leq i\leq p$, and 
consider the corresponding block $[{\bf s}, {\bf t}]=\prod_{i=1}^p [s_i,t_i]$. 
The dimension of $[{\bf s}, {\bf t}]$ is given by the number of indices $i$ such that $s_i<t_i$.
The generalized increment of a random field $X=(X_{\bf t})_{{\bf t}\in\real^p}$ on a block $[{\bf s}, {\bf t}]$ in $\real^p$ is defined (up to a factor $\pm 1$) by
\begin{equation}
\label{eq:XB}
X([{\bf s}, {\bf t}]):=\sum_{{\mathbf \epsilon}} (-1)^{p-\sum_{i=1}^p \epsilon_i}X(s_1+\epsilon_1(t_1-s_1),\cdots,s_p+\epsilon_p(t_p-s_p))
\end{equation}
where the sum above runs over ${\mathbf \epsilon}\in\{0,1\}^p$ such that $\epsilon_i=0$ if $s_i=t_i$.
Similarly for measure $\mu_{\bf t}$ with density $\phi_{{\bf t}}$, we define the increment on a block $[{\bf s}, {\bf t}]$ of $\phi_{{\bf t}}$ by
\begin{equation}
\label{eq:phiB}
\phi_{[{\bf s}, {\bf t}]}:=\sum_{{\mathbf \epsilon}} (-1)^{p-\sum_{i=1}^p\epsilon_i} 
\phi_{s_1+\epsilon_1(t_1-s_1),\dots,s_p+\epsilon_p(t_p-s_p)}
\end{equation}
where again the sum runs over  ${\mathbf \epsilon}\in\{0,1\}^p$ such that $\epsilon_i=0$ if $s_i=t_i$. 
Such generalized increments are easy to handle in our context and Example \ref{exa:2} below confirms, for uniform-type measures, that such increments make sense in our setting. 
In order to control such increments, we introduce our main condition on the densities $\phi_{{\bf t}}$ of $\mu_{{\bf t}}$: we say that the family of densities $(\phi_{\bf t})_{{\bf t}\in\real^p}$ satisfies property $(P_\gamma)$ for $\gamma\geq 1$ 
if for all $T>0$, there exists some constant $C_T>0$ such that for any $[{\bf s}, {\bf t}]\subset [-T,T]^p$, 
\begin{equation}
\tag{$P_\gamma$}
\|\phi_{[{\bf s}, {\bf t}]}\|_\gamma^\gamma \leq C_T \prod_{i:s_i<t_i} |t_i-s_i|.
\end{equation}
Here, $\|\phi\|_\gamma$ stands for the $L^\gamma(\real^p)$-norm of $\phi$.
The following examples justify that condition $(P_\gamma)$ is natural.  

\begin{example}\label{exa:1} Let $(\mu_{\bf t})_{{\bf t}\in\real^p}$ be the family of (signed) uniform measures  on the blocks $[{\bf 0},{\bf t}]=\prod_{i=1}^d[0,t_i]$, more precisely, $\phi_{\bf t}={\rm sign}(t_1)\cdots {\rm sign}(t_p)\ind_{[{\bf 0},{\bf t}]}$. 
Then, we verify that for any non-degenerated block $[s,t]\subset\real^p$, $\phi_{[{\bf s}, {\bf t}]}=\ind_{[{\bf s}, {\bf t}]}$ almost everywhere so that
$$
\|\phi_{[{\bf s}, {\bf t}]}\|_\gamma^\gamma= \prod_{1\leq i\leq p} |t_i-s_i|.
$$
In the case of a degenerated block $[s,t]\subset [-T,T]^p$, we have
$$
\|\phi_{[{\bf s}, {\bf t}]}\|_\gamma^\gamma= \prod_{i:s_i=t_i}|s_i|  \prod_{i:s_i<t_i} |t_i-s_i|\leq C_T \prod_{i:s_i<t_i} |t_i-s_i|
$$
with $C_T=\max(1,T)^p$. Hence $(P_\gamma)$ holds true for all $\gamma\geq 1$. 
Such uniform-type measures are used to analyze cumulative workload in one-dimensional model, see \cite{MRRS}, \cite{KT07}.  
\end{example}

\begin{example}
\label{exa:2}
 Let $(\mu_{\bf t})_{{\bf t}\in\real^p}$ be a family of measures with densities $\phi_{\bf t}$ 
such that for all $I\subset \{1,\dots, p\}$ and all $T>0$, 
\begin{equation}
\label{eq:ex2}
\Big\|\sup_{t\in [-T,T]^p}\partial_I \phi_{\bf t}({\bf y})\Big\|_\gamma<+\infty 
\end{equation}
where $\partial_I$ is the differential operator defined for $I=\{i_1,\dots,i_k\}$ by
$\partial_I =\partial ^k/\partial t_{i_1}\cdots\partial t_{i_k}$.
The following condition $(P'_\gamma)$ (that implies condition $(P_\gamma)$) is satisfied: 
for any $T>0$, there exists some constant $C_T>0$ such that for any $[{\bf s}, {\bf t}]\subset [-T,T]^p$, 
\begin{equation}
\tag{$P^{'}_\gamma$} 
\|\phi_{[{\bf s}, {\bf t}]}\|_\gamma^\gamma \leq C_T\prod_{i:s_i<t_i} |t_i-s_i|^\gamma.
\end{equation}
This fact is justified in Section \ref{sec:tightness}.
\end{example}

\medskip
\noindent
We now state our results for the functional convergences in the large ball and intermediate ball regime. 
We recall that the limits $Z_\alpha$ and $J$ below are defined in Proposition~\ref{prop:fdd}. 
\begin{theo}
\label{theo:CV_functional1}
Suppose conditions $(\bf A)$ hold. 
Let $\mu_{\bf t}(d{\bf y})=\phi_{\bf t}({\bf y})d{\bf y}, {\bf t}\in\real^p$, be a parametric family of measures in $L^1(\real^d)\cap L^\alpha(\real^d)$ satisfying conditions $(P_1)$ and $(P_\alpha)$.
\begin{enumerate}
\item (Large ball regime) If $\lambda(\rho) \rho^\beta\to+\infty$, then, setting $n(\rho)=(\lambda(\rho) \rho^\beta)^{1/\alpha}$, we have as $\rho\to 0$: 
\begin{equation*}
\label{eq:stable}
\widetilde M_\rho(\mu_{\bf t}) \stackrel{{\mathcal C}(\real^p)}{=\!\!=\!\!\Longrightarrow} Z_\alpha(\mu_{\bf t}), \quad {\bf t}\in\real^p.
\end{equation*}
\item (Intermediate ball regime) If $\lambda(\rho)\rho^\beta \to a>0$, then, setting $n(\rho)=1$, we have as $\rho\to 0$:
$$
\widetilde M_\rho(\mu_{\bf t})  \stackrel{{\mathcal C}(\real^p)}{=\!\!=\!\!\Longrightarrow}J_a(\mu_{\bf t}), \quad {\bf t}\in\real^p.
$$
\end{enumerate}
\end{theo}

As a by-product of the moment estimates used to prove tightness, we obtain H\"older-regularity properties in the case $\alpha=2$. This is the content of the following result. 
\begin{prop}
\label{corol:CV_functional3}
Suppose conditions $(\bf A)$ hold with $\alpha=2$. 
\begin{enumerate}
\item If $G$ has a finite variance 
and the family of measures $\mu_{\bf t}(d{\bf y})=\phi_{\bf t}({\bf y})d{\bf y}$ satisfies  $(P_1)$ and $(P_2)$, then the Gaussian limit process $(Z_2(\mu_{\bf t}))_{{\bf t}\in\real^p}$ of Theorem \ref{theo:CV_functional1} is $\gamma$-H\"older for all  $\gamma<\frac{3d-\beta}{2d}$.
\item If there is $k\geq 2p$ such that  $h\in L^k(\real^d)$, $G$ has finite moment of order $k$, 
and the family of measures $\mu_{\bf t}(d{\bf y})=\phi_{\bf t}({\bf y})d{\bf y}$ satisfies $(P_1)$ and $(P_k)$, 
then the limit process $(J_a(\mu_{\bf t}))_{{\bf t}\in\real^p}$ of Theorem \ref{theo:CV_functional1} is $\gamma$-H\"older for all $\gamma<\frac{3d-\beta}{2d}-\frac pk$. 
\end{enumerate}
\end{prop}
 
Observe that in dimension $d=1$, we recover at the limit the fractional Brownian motion obtained at the limit in \cite{KT07} 
with the H\"older-regularity $\gamma<(3-\beta)/2\in (\frac 12,1)$. 

\medskip
\noindent
Finally, we consider functional convergence in the space of distribution. 
Let $\cD(\real^d)$ be the space of smooth compactly supported functions, and $\cD'(\real^d)$ be its dual, {\it i.e.} the space of distributions. 
We show that, for all $\rho>0$, $\widetilde M(\rho)$ can be seen as a random distribution and state functional convergence in $\cD'(\real^d)$.
\begin{theo}
\label{theo:distrib} 
Suppose conditions $(\bf A)$ hold. 
\begin{enumerate}
\item For each $\rho>0$, $M_\rho$ induces a random distribution, {\it i.e.} the linear form
$$
M_\rho:\left\{\begin{array}{ccl} \cD(\real^d)&\to&\real \\ \phi&\mapsto& M_\rho(\phi({\bf y})d{\bf y})\end{array}\right.
$$
is almost surely continuous.
\item (Large ball regime) If $\lambda(\rho) \rho^\beta\to+\infty$, then, setting $n(\rho)=(\lambda(\rho) \rho^\beta)^{1/\alpha}$ in $\widetilde M_\rho(\mu)$, we have as $\rho\to 0$: 
$$
\widetilde M_\rho(\mu) \stackrel{{\mathcal D}'(\real^p)}{=\!\!=\!\!\Longrightarrow} Z_\alpha(\mu),\quad \mu\in\cD(\real^d).
$$
\item (Intermediate ball regime) If $\lambda(\rho)\rho^\beta \to a>0$, then, setting $n(\rho)=1$ in $\widetilde M_\rho(\mu)$, we have: 
$$
\widetilde M_\rho(\mu)\stackrel{{\mathcal D}'(\real^p)}{=\!\!=\!\!\Longrightarrow} J_a(\mu), \quad \mu\in\cD(\real^d).
$$
\end{enumerate}
\end{theo}
\begin{Rem}
{\rm 
The result is still true if $\cD(\real^d)$ is replaced by ${\mathcal C}^k_K(\real^d)$ the space of compactly supported functions of class ${\mathcal C}^k$ on $\real^d$ 
and $\cD'(\real^d)$ is replaced by the dual of ${\mathcal C}^k_K(\real^d)$, $k\in\nit\setminus\{0\}$.
}
\end{Rem}


\section{Proofs}
\label{sec:proof}
The proof of Theorem \ref{theo:CV_functional1}, as usual for functional convergences, consists of two arguments: 
{\it fdd} convergences and tightness. 
The first one is given in Proposition \ref{prop:fdd} ; this is a slight modification of \cite{BD} and the changes are discussed in Section~\ref{sec:fdd}.
The proof of tightness is given in Section~\ref{sec:tightness} and it relies on moment estimates previously obtained in Section \ref{sec:moment}. 
H\"older regularity also relies on moments and cumulants estimates and Proposition \ref{corol:CV_functional3} is proved in Section \ref{sec:regularity}. 
Section \ref{sec:distribution} is devoted to the proof of functional convergence in the space of distributions $\cD'(\real^d)$. 

Recall that throughout the paper, we assume that conditions $({\bf A})$ hold, and in particular we consider $d<\beta<\alpha d$. All the asymptotics are considered as  $\rho\to~0$. 

\subsection{{\it fdd} convergences}
\label{sec:fdd}

The results of \cite{BD} do not apply directly since the model we consider here is not exactly the same: in \cite{BD}, the probability $F$ is assumed to have a density, the fading function $h$ is replaced by $\ind_{B({\bf 0},1)}$ and the {\it fdd} convergence is proved on the space ${\mathcal M}_{\alpha, \beta}$ of measure $\mu$ satisfying \eqref{eq:Malphabeta}. 
Here the radius distribution $F$ and the shape function $h$ are more general (see conditions (${\bf A}$)), but we obtain convergence only on the space of measures $\mu(d{\bf y})=\phi({\bf y})d{\bf y}$ with $\phi\in L^1(\real^d)\cap L^\alpha(\real^d)$. 
However, as briefly explained below, {\it verbatim} changes in the proofs of \cite{BD} shows that the results of \cite{BD} still apply in our context. 

Indeed, first a careful reading of the proofs in \cite{BD} shows that the existence of the density $f$ of $F$ is used only  in Lemma 3.2 therein. 
But this lemma can be replaced by Lemmas~2 and 3 in \cite{KLNS} deriving the same result as Lemma 3.2 in \cite{BD} under the weaker assumption \eqref{eq:F2}. 
Observe in particular that $\bar F_\rho(1)=\bar F(1/\rho)\sim C_\beta \rho^\beta$ and that the continuity requirement in Lemmas~2 and 3 in \cite{KLNS} is ensured in our setting by condition \eqref{eq:halpha} and Lemma \ref{lem:cont} below.  

Second, the bound \eqref{eq:Malphabeta} can be replaced by the following condition \eqref{eq:remplac}: 
for $h\in L^1(\real^d)\cap L^\alpha(\real^d)$ and $\mu(d{\bf y})=\phi({\bf y})d{\bf y}$ with $\phi\in L^1(\real^d)\cap L^\alpha(\real^d)$, we have
\begin{equation}
\label{eq:remplac}
\int_{\real^d} |\mu[\tau_{{\bf x},r}h]|^\alpha d{\bf x}\leq C (r^d\wedge r^{\alpha d}),
\end{equation}
this is justifies in Lemma \ref{lemme:muh} in the Appendix with $\gamma:=\alpha$ therein.

\medskip
As a consequence, Theorems 2.4 and 2.11 in \cite{BD} remain true in our setting and ensure the {\it fdd} convergences stated in Proposition \ref{prop:fdd}.  
 

\subsection{Moment estimates}
\label{sec:moment}
As we will see, our results on tightness and H\"older regularity strongly rely on moment estimates for the rescaled random field $\widetilde M_\rho$. 
Since the following properties of the moment are also interesting in their own right, 
they are stated in the following proposition.
\begin{prop}
\label{prop:moment} 
Suppose conditions $(\bf A)$ hold. 
\begin{enumerate}
\item Let $0<\gamma<\alpha$. There exists some constant $C:=C(F,G,h,\alpha,\beta,\gamma,d)$ depending only on $F,G,h,\alpha,\beta,\gamma,d$ ({\it i.e.} not on $\rho$ and $\phi$) 
such that for all $\mu(d{\bf y})=\phi({\bf y})d{\bf y}\in L^1(\real^d)\cap L^\alpha(\real^d)$ and all $\rho>0$, we have
\begin{equation}
\label{eq:moment1a}
\ee\left[\left|\widetilde M_\rho(\mu)\right|^\gamma\right]
\leq C \left[\frac{\lambda(\rho)\rho^\beta}{n(\rho)^\alpha}\right]^{\gamma/\alpha}
\|\phi\|_\alpha^{\frac{\gamma(\beta-d)}{(\alpha-1)d}}\|\phi\|_1^{\frac{\gamma(\alpha d-\beta)}{(\alpha-1)d}}.
\end{equation}
\item Suppose that $\alpha=2$ and $G$ has a finite moment of order $k\in \nit\setminus\{0,1\}$, $h\in L^1(\real^d)\cap L^k(\real^d)$. 
There exists some constant $C:=C(F,G,h,\alpha,\beta,k,d)$ depending only on $F,G,h,\alpha,\beta,k,d$ 
such that for all $\mu(d{\bf y})=\phi({\bf y})d{\bf y}\in L^1(\real^d)\cap L^k(\real^d)$ and all $\rho>0$,
\begin{equation}
\label{eq:cumulant1}
\left|c_k\big(\widetilde M_\rho(\mu)\big)\right|
\leq C \frac{\lambda(\rho)\rho^\beta}{{n(\rho)^k}}
\|\phi\|_k^{\frac{k(\beta-d)}{(k-1)d}}\|\phi\|_1^{\frac{k(k d-\beta)}{(k-1)d}},
\end{equation}
where $c_k\big(\widetilde M_\rho(\mu)\big)$ is the cumulant of order $k$ of $\widetilde M_\rho(\mu)$. 
\end{enumerate}
\end{prop}

As a by-product of these moment estimates and the finite-dimensional convergence (Proposition \ref{prop:fdd}), we obtain the following result stating the convergence of moments:
\begin{corol}
\label{corol:moment}
Suppose conditions ({\bf A}) hold and $\mu(d{\bf y})=\phi({\bf y})d{\bf y}$ with $\phi\in L^1(\real^d)\cap L^\alpha(\real^d)$.
Let $0<\gamma<\alpha$.
\begin{enumerate}
\item If $\lambda(\rho)\rho^\beta\to +\infty$, then, setting $n(\rho)=(\lambda(\rho)\rho^\beta)^{1/\alpha}$ in $M_\rho(\mu)$, we have as $\rho\to 0$:
$$
\ee\left[\left|\widetilde M_\rho(\mu)\right|^\gamma\right]\to \ee\left[\left|Z_\alpha(\mu)\right|^\gamma\right].
$$
\item If $\lambda(\rho)\rho^\beta\to a$, then, setting $n(\rho)=1$ in $M_\rho(\mu)$, we have as $\rho\to 0$:
$$
\ee\left[\left|\widetilde M_\rho(\mu)\right|^\gamma\right]\to \ee\left[\left|J_a(\mu)\right|^\gamma\right].
$$
\end{enumerate}
In the case $\alpha=2$, if furthermore $G$ has a finite moment of order $n\geq 2$ and $h,\phi\in L^1(\real^d)\cap L^n(\real^d)$, then the above convergence of moments holds for all $0\leq \gamma\leq n$.
\end{corol}


\medskip
\noindent
{\bf Proof of Proposition \ref{prop:moment}.}\\
{\bf First point:} The estimate for $\ee\left[\left|\widetilde M_\rho(\mu)\right|^\gamma\right]$ relies on the following expression of the fractional moment 
(see \cite{vBE} or \cite[Eq. (60)]{KT07}): 
if $X$ is a random variable with characteristic function $\varphi_X(t)=\ee[\exp(itX)]$, then we have, for $1<\gamma<2$, 
\begin{equation}
\label{eq:mgamma}
\ee[|X|^\gamma]=A(\gamma)\int_0^{+\infty}(1-|\varphi_X(\theta)|^2)\theta^{-1-\gamma}d\theta
\end{equation}
where
$$
A(\gamma)=\left(\int_0^{+\infty}(1-\cos(x))x^{-1-\gamma}d{\bf x}\right)^{-1}<+\infty.
$$
Since $\widetilde M_\rho(\mu)$ is a Poisson integral, its caracteristic function is given by (see Lemma \ref{prop:cumulant})
\begin{equation*}
\varphi_{\widetilde M_\rho(\mu)}(\theta)
=\exp\left(\int_{\real^d\times\real^+} \Psi_G\left(n(\rho)^{-1}\theta \mu[\tau_{{\bf x},r}h ]\right)\lambda(\rho)d{\bf x}F_\rho(dr)\right)
\end{equation*}
with $\Psi_G(u)=\int_\real (e^{ium}-1-ium )G(dm)$. Hence,
\begin{eqnarray}
&&1-|\varphi_{\widetilde M_\rho(\mu)}(\theta)|^2 \nonumber\\
&=&1-\exp\left(\int_{\real^d\times\real^+} 2{\rm Re}\Big(\Psi_G\left(n(\rho)^{-1}\theta \mu[\tau_{{\bf x},r}h ]\right)\Big)\lambda(\rho)d{\bf x}F_\rho(dr)\right) \nonumber\\
&\leq&1-\exp\left(-2\int_{\real^d\times\real^+} \left|\Psi_G\left(n(\rho)^{-1}\theta \mu[\tau_{{\bf x},r}h ]\right)\right|\lambda(\rho)d{\bf x}F_\rho(dr)\right)\nonumber \\
&\leq &1-\exp\left(-2C(G)n(\rho)^{-\alpha}\lambda(\rho)|\theta|^\alpha\int_{\real^d\times\real^+} \left|\mu[\tau_{{\bf x},r}h ]\right|^\alpha d{\bf x}F_\rho(dr)\right)
\label{eq:bound1}.
\end{eqnarray}
Using Lemma \ref{lemme:intmuh} in the Appendix 
with $\gamma:=\alpha>\beta/d$, we have
\begin{equation}
\label{eq:bound2}
\int_{\real^d\times\real^+} |\mu[\tau_{{\bf x},r}h]|^\alpha d{\bf x}F_\rho(dr)\leq 
\rho^\beta \widetilde C(\phi)
\end{equation}
where $\widetilde C(\phi)$ is given by
\begin{eqnarray*}
\widetilde C(\phi)&=&C(F)\frac{(\alpha-1)\beta d}{(\alpha d-\beta)(\beta-d)}
(\|\phi\|_\alpha \|h\|_1)^{\frac{\alpha(\beta-d)}{(\alpha-1)d}}
(\|\phi\|_1\|h\|_\alpha)^{\frac{(\alpha d-\beta)\alpha}{(\alpha-1)d}}
\end{eqnarray*}
and $C(F)$ is a constant depending only on $F$. 
Plugging the bounds \eqref{eq:bound1} and \eqref{eq:bound2} in \eqref{eq:mgamma}, we obtain
\begin{eqnarray}
\nonumber
\ee\left[\left|\widetilde M_\rho(\mu)\right|^\gamma\right]
&\leq& A(\gamma)\int_0^{+\infty} \left(1-\exp\left(-2C(G)n(\rho)^{-\alpha}\lambda(\rho)\rho^\beta |\theta|^\alpha \tilde C(\phi) \right)\right) \theta^{-1-\gamma}d\theta\\
\label{eq:boundEM}
&=&A(\gamma)A(\alpha,\gamma)\left(\frac{\lambda(\rho)\rho^\beta}{n(\rho)^{\alpha}}\right)^{\gamma/\alpha}(2C(G)\widetilde C(\phi))^{\gamma/\alpha}
\end{eqnarray}
with a straightforward change of variables in \eqref{eq:boundEM} and 
$$
A(\alpha,\gamma)=\int_0^{+\infty} \left(1-\exp\left(-\theta^\alpha  \right)\right) \theta^{-1-\gamma}d\theta<+\infty.
$$
This gives the result \eqref{eq:moment1a} with the constant
\begin{eqnarray*}
&&C(F,G,h,\alpha,\beta,\gamma,d)\\
&=&A(\gamma)A(\alpha,\gamma)\left(2C(G) C(F)\frac{(\alpha-1)\beta d}{(\alpha d-\beta)(\beta-d)}
 \|h\|_1^{\frac{\alpha(\beta-d)}{(\alpha-1)d}}
\|h\|_\alpha^{\frac{(\alpha d-\beta)\alpha}{(\alpha-1)d}}\right)^{\gamma/\alpha}.
\end{eqnarray*}
The case $1<\gamma<\alpha$ is proved. 
The case $0<\gamma\leq 1$ comes from the previous case applied to any $\gamma'\in(1,\alpha)$ combined with the Jensen inequality which implies 
$\ee[|X|^\gamma]\leq \ee[|X|^{\gamma'}]^{\gamma/\gamma'}$. 

\medskip
\noindent
{\bf Second point:} 
We use the general form of the cumulant of a Poisson integral as recalled in Lemma \ref{prop:cumulant}. 
Using Lemma \ref{lemme:intmuh} with $\gamma:=k \geq 2>\beta/d$, we have:
\begin{eqnarray*}
&&c_k(\widetilde M_\rho(\mu))\\
&=&\frac{\lambda(\rho)}{{n(\rho)}^{k}}\left(\int_{\real}m^kG(dm)\right) \left(\int_{\real^d\times\real^+} [\mu[\tau_{{\bf x},r}h]]^k d{\bf x}F_\rho(dr)\right)\\
&\leq&\!\!\frac{\lambda(\rho)\rho^\beta}{{n(\rho)^k}} C(F)  \frac{(k-1)\beta d}{(k d-\beta)(\beta-d)}
(\|\phi\|_k \|h\|_1)^{\frac{k(\beta-d)}{(k-1)d}}(\|\phi\|_1\|h\|_k)^{\frac{(k d-\beta)k}{(k-1)d}}\left(\int_{\real}\!\! m^kG(dm)\right).
\end{eqnarray*}
This gives the bound \eqref{eq:cumulant1} with the constant  
\begin{eqnarray*}
&&C(F,G,h,\alpha,\beta,k,d)\\
&=&C(F) \frac{(k-1)\beta d}{(k d-\beta)(\beta-d)}
\|h\|_1^{\frac{k(\beta-d)}{(k-1)d}}\|h\|_k^{\frac{(k d-\beta)k}{(k-1)d}}\left(\int_{\real}m^kG(dm)\right).
\end{eqnarray*}
\CQFD
\ \\
\noindent
{\bf Proof of Corollary \ref{corol:moment}.}\\
We use the following basic result: if $X_n$ weakly converges to $X$ as $n\to+\infty$ and $\ee\left[|X_n|^{\gamma'}\right]$ is bounded, 
then the convergence of moments  $\ee\left[|X_n|^{\gamma}\right]\to\ee\left[|X|^{\gamma}\right]$ holds for any $0<\gamma<\gamma'$ 
(in this case, the family $|X_n|^{\gamma}$ is indeed bounded in $L^{\gamma'/\gamma}$ and hence equi-integrable). 
As a consequence of this result, the finite-dimensional convergence stated in Proposition \ref{prop:fdd} and the moment estimates obtained in Proposition \ref{prop:moment} yield the convergence of moments for $0<\gamma<\alpha$ and $0<\gamma<n$ respectively. 
The case $\gamma=n$ should be proved separately but we omit the details.
\CQFD 


\subsection{Tightness}
\label{sec:tightness}

This section is devoted to the second step in the proof of Theorem \ref{theo:CV_functional1}, {\it i.e.} tightness. 
More precisely, we show that:
\begin{prop}
\label{prop:tightness}
Under the assumptions of Theorem \ref{theo:CV_functional1}, 
the family of random fields $(\widetilde M_\rho(\mu_{\bf t}))_{{\bf t}\in\real^p}, \rho\leq 1$, is tight in ${\mathcal C}(\real^p)$. 
\end{prop}

The proof of Proposition \ref{prop:tightness} relies on a suitable control of the moment of the generalized increments of $\widetilde M_\rho(\mu_{\bf t})$
with the following Censov criterion for which we refer to \cite[p. 16]{DP} or to \cite{BW}.
\begin{prop}
\label{prop:tightness_DP}
Let $(X_n)_{n\in\nit}$ be a sequence of random fields on $\real^p$ such that:  
\begin{enumerate}
\item The family of random variable $(X_n({\bf 0}))_{n\in\nit}$ is tight; 
\item For all $T>0$, there are constants $\gamma> 0$, $\delta>1$ and $C_T>0$ such that, for all $[{\bf s}, {\bf t}]\subset [-T,T]^p$ and $n\in\nit$,
\begin{eqnarray*}
\label{eq:tight_DP0}
\ee\left[|X_n([{\bf s}, {\bf t}])|^\gamma\right]\leq C_T\prod_{i:s_i\neq t_i}|t_i-s_i|^\delta. 
\end{eqnarray*}
\end{enumerate}
\end{prop}
\begin{Rem}
{\rm Proposition \ref{prop:tightness_DP} is a weaker version of the original criterion given in \cite{DP} that we have adapted to our setting. 
Note that $\prod_{i:s_i\neq t_i}|t_i-s_i|$ is the $k$-dimensional Lebesgue measure of $[{\bf s}, {\bf t}]$, 
where $k$ is the true dimension of the block $[{\bf s}, {\bf t}]$. 
In \cite{DP}, the Lebesgue measure is replaced by general Radon measure with diffuse marginals. 
Originally, the criterion is expressed in terms of bounds on the tails $\bbP(|X_n([{\bf s}, {\bf t}])|>x)$ 
instead of the moments $\ee\left[|X_n([{\bf s}, {\bf t}])|^\gamma\right]$. 
Moreover, it is explained in \cite{DP} that when $k$-dimensional faces are involved, it is enough to check the condition on blocks $[{\bf s}, {\bf t}]$ such that $s_i=t_i=0$ for the degenerated dimensions. 
For the sake of clearness, we do not insist on such general conditions.
}
\end{Rem}

\medskip
\noindent
{\bf Proof of Proposition \ref{prop:tightness}.}\\
We apply Proposition \ref{prop:tightness_DP} to the family $X_\rho({\bf t})=\widetilde M_\rho(\mu_{\bf t})$, $\rho\leq 1$. 
To that aim, observe first that the tightness of $X_\rho({\bf 0})$ is a consequence of the one-dimensional convergence stated in Proposition \ref{prop:fdd}. 
Furthermore, using the definitions \eqref{eq:XB} and \eqref{eq:phiB} together with the linearity of the generalized random field $\widetilde M_\rho$, we see easily that 
$$
X_\rho([{\bf s}, {\bf t}])=\widetilde M_\rho(\mu_{[{\bf s}, {\bf t}]})\quad {\rm with}\quad \mu_{[{\bf s}, {\bf t}]}(d{\bf y})=\phi_{[{\bf s}, {\bf t}]}({\bf y})d{\bf y}.
$$
Proposition \ref{prop:moment}$-1)$ gives the following moment estimate for $0<\gamma <\alpha$:
\begin{eqnarray}
\nonumber
\ee\left[|X_\rho([{\bf s}, {\bf t}])|^\gamma\right]&=&\ee\left[|\widetilde M_\rho(\mu_{[{\bf s}, {\bf t}]})|^\gamma\right]\\
\label{eq:increments1}
&\leq &  C \left[\frac{\lambda(\rho)\rho^\beta}{n(\rho)^\alpha}\right]^{\gamma/\alpha}
\|\phi_{[{\bf s}, {\bf t}]}\|_\alpha^{\frac{\gamma(\beta-d)}{(\alpha-1)d}}\|\phi_{[{\bf s}, {\bf t}]}\|_1^{\frac{\gamma(\alpha d-\beta)}{(\alpha-1)d}}.
\end{eqnarray}
Using $(P_1)$ and $(P_\alpha)$ for $\phi_{[{\bf s}, {\bf t}]}$, we obtain:
\begin{eqnarray}
\nonumber
\ee\left[|X_\rho([{\bf s}, {\bf t}])|^\gamma\right]&\leq& C \left[\frac{\lambda(\rho)\rho^\beta}{n(\rho)^\alpha}\right]^{\gamma/\alpha}
\left(\prod_{i:s_i< t_i}|t_i-s_i|\right)^{\frac{\gamma(\beta-d)}{\alpha(\alpha-1)d}+\frac{\gamma(\alpha d-\beta)}{(\alpha-1)d}}\\
\label{eq:increments2}
&\leq & C\left[\frac{\lambda(\rho)\rho^\beta}{n(\rho)^\alpha}\right]^{\gamma/\alpha}
\left(\prod_{i:s_i< t_i}|t_i-s_i|\right)^{\frac{\gamma}{\alpha}(1+\alpha-\beta/d)}.
\end{eqnarray}
Now, observe that the sequence $\lambda(\rho)\rho^\beta/n(\rho)^\alpha$ is bounded since under the two asymptotics investigated 
$\lambda(\rho)\rho^\beta/n(\rho)^\alpha$ converge to some finite constant as $\rho\to 0$. 
Furthermore since $d<\beta<\alpha d$ and $1+\alpha-\beta/d>1$, the exponent $\frac{\gamma}{\alpha}(1+\alpha-\beta/d)$ is (strictly) larger than $1$ for $\gamma$ close enough to $\alpha$. 
This proves the second condition in Proposition~\ref{prop:tightness_DP}. 
Proposition \ref{prop:tightness} and thus Theorem \ref{theo:CV_functional1} easily follow.
\CQFD 

\medskip
We finish this section with the proof of Example \ref{exa:2} where the differentiability condition \eqref{eq:ex2} 
is stated to be sufficient for $(P_\gamma')$. 

\medskip
\noindent
{\bf Proof for Example \ref{exa:2}.}
Let ${\bf s}\leq {\bf t}$. 
We first assume that, for all $1\leq i\leq p$, $ -T\leq s_i<t_i\leq T$. 
We have
$$
\phi_{[{\bf s}, {\bf t}]}({\bf y})=\int_{[{\bf s}, {\bf t}]} \partial_{\{1,\dots,p\}}\phi_{\bf u}({\bf y})d{\bf u}. 
$$
Using H\"older inequality, we have:
\begin{eqnarray*}
\left|\int_{[{\bf s}, {\bf t}]}\partial_{\{1,\dots,p\}}\phi_{\bf u}({\bf y})d{\bf u}\right|^\gamma
&\leq& |[{\bf s}, {\bf t}]|^{\gamma-1} \int_{[{\bf s}, {\bf t}]} \left|\partial_{\{1,\dots,p\}}\phi_{\bf u}({\bf y})\right|^\gamma d{\bf u}\\
&\leq& |[{\bf s}, {\bf t}]|^{\gamma}\sup_{t\in [-T,T]^p}|\partial_{\{1,\dots,p\}} \phi_{\bf u}({\bf y})|^\gamma
\end{eqnarray*}
so that
$$
\|\phi_{[{\bf s}, {\bf t}]}\|_\gamma^\gamma 
\leq \Big\|\sup_{t\in [-T,T]^p}\partial_{\{1,\dots,p\}} \phi_{\bf t}({\bf y})\Big\|_\gamma^\gamma \ \prod_{i=1}^p |t_i-s_i|^\gamma.
$$
In general, if ${\bf s}\leq {\bf t}$, let $I$ be the set of indices such that $s_i<t_i$. 
We show similarly that 
$$
\|\phi_{[{\bf s}, {\bf t}]}\|_\gamma^\gamma 
\leq \Big\|\sup_{t\in [-T,T]^p}\partial_{I} \phi_{\bf t}({\bf y})\Big\|_\gamma^\gamma \ \prod_{i\in I} |t_i-s_i|^\gamma.
$$
\CQFD


\subsection{H\"older-regularity when $\alpha=2$} 
\label{sec:regularity}
In this section, we prove Proposition \ref{corol:CV_functional3} where H\"older-regularity is stated for the limit $Z_2(\mu_{\bf t})$ (resp.  $J_a(\mu_{\bf t})$) when $G$ has finite variance (resp. finite moments of any order). 
In particular, we set $\alpha=2$ and we assume $\beta\in(d,2d)$.
H\"older-regularity is proven using again moment estimates for increments. 
However since we have not found in the literature H\"older-regularity result relying on generalized increments in dimension $p\geq 1$, 
we use standard increments $Z_2(\mu_{\bf t})-Z_2(\mu_{\bf s})$ (resp. $J_a(\mu_{\bf t})-J_a(\mu_{\bf s})$) 
that will be controled thank to the following elementary observation: suppose the family $(\phi_{\bf t})_{{\bf t}\in\real^p}$ satisfies condition $(P_\gamma)$, then for all $T>0$, for all $[{\bf s},{\bf t}]\subset[-T,T]^p$,
\begin{equation}
\label{eq:Pmod} 
\|\phi_{\bf s}-\phi_{\bf t}\|_\gamma^\gamma \leq C_T p^\gamma \|{\bf s}-{\bf t} \|_{\infty}^\gamma.
\end{equation}
To see this, rewrite a standard increment as a sum of $d$ block increments along axis.
More precisely, for ${\bf s}, {\bf t}\in\real^p$, define ${\bf u}^i\in \real^p$ by $u^i_j=t_j$ for $j\leq i$ and $u^i_j=s_j$ for $j>i$ so that, for each $0\leq i\leq p-1$, ${\bf u}^i$ and ${\bf u}^{i+1}$ only differ from at most one coordinate. 
For $\gamma>1$, we have
$$
\|\phi_{\bf s}-\phi_{\bf t}\|_\gamma^\gamma
\leq p^{\gamma-1} \sum_{i=1}^p \big\|\phi_{[{\bf u}^{i-1},{\bf u}^{i}]}\big\|_\gamma^\gamma
$$
so that condition $(P_\gamma)$ implies
$$
\|\phi_{\bf s}-\phi_{\bf t}\|_\gamma^\gamma
\leq C_Tp^{\gamma-1}\  \sum_{i=1}^p |t_i-s_i|
\leq C_Tp^\gamma \ \|{\bf s}-{\bf t}\|_\infty. 
$$
\begin{Rem}
{\rm The trick consisting in controlling the moment of standard increments by generalized increments in \eqref{eq:Pmod} cannot be used in Section \ref{sec:tightness}. 
Indeed, similar bounds as in \eqref{eq:increments1} but for standard increments $\widetilde M_\rho(\mu_{\bf t})-\widetilde M_\rho(\mu_{\bf s})$ 
combined with \eqref{eq:Pmod} would yield
$$
\ee[|\widetilde M_\rho(\mu_{\bf t})-\widetilde M_\rho(\mu_{\bf s})|^\gamma]
\leq CC_Tp^\gamma\left[\frac{\lambda(\rho)\rho^\beta}{n(\rho)^\alpha}\right]^{\gamma/\alpha}
\|{\bf s}-{\bf t}\|_\infty^{\frac{\gamma}{\alpha}(1+\alpha-\beta/d)}. 
$$
instead of \eqref{eq:increments2} which is far from a bound in $\|{\bf s}-{\bf t}\|^\delta$ with $\delta>p$ required in the criterion for tightness relying on standard increments (see \cite[Th. 1.4.1]{Kun}). 
The trick works when $\alpha=2$ because in this case, the limit is Gaussian and we can artificially increase the exponent in the bound thank to the 
relation between moments and variance, see~\eqref{eq:momentgaussien}. 
}
\end{Rem}
\medskip
\noindent
{\bf Proof of Proposition \ref{corol:CV_functional3}.}\\
{\bf First point:} We assume that $G$ has a finite variance, $h\in L^1(\real^d)\cap L^2(\real^d)$ and $(P_\gamma)$ holds for $\gamma=1,2$. 
Using Lemma \ref{lemme:muh}, we have, for $\mu(d{\bf y})=\phi({\bf y})d{\bf y}$,
\begin{eqnarray*}
\Var(Z_2(\mu))&=&\int_{\real^d} \left|\mu[\tau_{{\bf x},r}h]\right|^2 \frac{\sigma^\alpha C_2}{r^{1+\beta}}drd{\bf x}\\
&\leq& \int_{\real^d} (r^d\|\phi\|_1^2  \|h\|_2^2 \wedge  r^{2d}\|\phi\|_2^2\|h\|_1^2)\sigma^\alpha C_2r^{-1-\beta}dr\\
&=& \sigma^\alpha C_2\left( \frac{\|\phi\|_2^2\|h\|_1^2}{2d-\beta}\left(\frac{\|\phi\|_1^2  \|h\|_2^2}{\|\phi\|_2^2\|h\|_1^2}\right)^{2-\beta/d}
+\frac{\|\phi\|_1^2  \|h\|_2^2}{\beta-d}\left(\frac{\|\phi\|_1^2  \|h\|_2^2}{\|\phi\|_2^2\|h\|_1^2}\right)^{1-\beta/d}\right)\\
&=& \frac{\sigma^\alpha C_2d}{(2d-\beta)(\beta-d)}\left(\|\phi\|_2^{\beta/d-1}\|h\|_1^{\beta/d-1}\|\phi\|_1^{2-\beta/d}\|h\|_2^{2-\beta/d}\right)^2\\
&\leq& C \|\phi\|_2^{2(\beta/d-1)}\|\phi\|_1^{2(2-\beta/d)}
\end{eqnarray*}
for some finite constant $C$. 
Next, the properties $(P_1)$ and $(P_2)$ together with \eqref{eq:Pmod} entail that for all ${\bf s}, {\bf t} \in\real^p$:
\begin{eqnarray*}
\Var(Z_2(\mu_{\bf s}-\mu_{\bf t}))
&\leq& C \|{\bf s}-{\bf t}\|_\infty^{3-\beta/d}.
\end{eqnarray*}
We derive easily the moment of order $2n$ of the Gaussian random variable $Z_2(\mu_{\bf s}-\mu_{\bf t})$,  
\begin{eqnarray}
\label{eq:momentgaussien}
\ee[Z_2(\mu_{\bf s}-\mu_{\bf t})^{2n}]&=&\frac{(2n-1)!}{(n-1)!2^{n-1}}\Var(Z_2(\mu_{\bf s}-\mu_{\bf t}))^n\\
\nonumber
&\leq& C \frac{(2n-1)!}{(n-1)!2^{n-1}} \|{\bf s}-{\bf t}\|_\infty^{(3-\beta/d)n}.
\end{eqnarray}
Using a standard regularity criterion for random fields (see e.g. \cite[Th. 1.4.1]{Kun}), $Z_2(\mu_{\bf t})$ has, almost surely, H\"older-continuous path for all index striclty less than $\frac{3d-\beta-pd/n}{2d}$.
Letting $n$ go to $+\infty$, we obtain that  $Z_2(\mu_{\bf t})$ is, almost surely, H\"older-continuous for all index striclty less than $\frac{3d-\beta}{2d}$.

\medskip
\noindent
{\bf Second point:} 
With $\mu=\mu_{\bf s}-\mu_{\bf t}$, Proposition \ref{prop:moment} entails 
\begin{equation}
\label{eq:cumulant2b}
\left|c_k(\widetilde M_\rho(\mu_{\bf s}-\mu_{\bf t}))\right|
\leq C \frac{\lambda(\rho)\rho^\beta}{{n(\rho)^k}}
\|\phi_{\bf s}-\phi_{\bf t}\|_k^{\frac{k(\beta-d)}{(k-1)d}}\|\phi_{\bf s}-\phi_{\bf t}\|_1^{\frac{k(k d-\beta)}{(k-1)d}},
\end{equation}
Now observe that in the intermediate regime, 
$\lambda(\rho)\rho^\beta/n(\rho)^k$ remains bounded since it converges to $a$ as $\rho\to 0$. 
Furthermore, using properties $(P_1)$ and $(P_k)$ together with \eqref{eq:Pmod} in \eqref{eq:cumulant2b}, we deduce
$$
\left|c_k(\widetilde M_\rho(\mu_{\bf s}-\mu_{\bf t}))\right|
\leq C \|{\bf s}-{\bf t}\|_\infty^{1+k-\beta/d}.
$$
Next, recall that the moments of a random variable $X$ are expressed in terms of its cumulants by the so-called complete Bell polynomials, 
{\it i.e.} $\ee[X^k]=B_k(c_1(X),\dots, c_k(X))$ with 
$$
B_k(c_1,\dots, c_k)=\sum_{i_1+2i_2\cdots+ki_k=k} K_n(i_1,\dots,i_n) c_1^{i_1}\dots c_k^{i_k}
$$ 
where $i_1,\dots, i_k$ are non-negative integers 
and $ K_k(i_1,\dots,i_k)$ are coefficients whose explicit (involved) form is not required in our argument. 
Since we have 
$$c_{1}(\widetilde M_\rho(\mu_{\bf s}-\mu_{\bf t}))^{i_1}\dots c_k(\widetilde M_\rho(\mu_{\bf s}-\mu_{\bf t}))^{i_k}=0$$
 when $i_1\not=0$, 
we can assume, without loss of generality, that $i_1=0$. 
For all ${\bf s},{\bf t}\in\real^d$ such that $\|{\bf s}-{\bf t}\|_\infty\leq 1$, we have 
\begin{eqnarray*}
\nonumber 
c_2(\widetilde M_\rho(\mu_{\bf s}-\mu_{\bf t}))^{i_2}\dots c_k(\widetilde M_\rho(\mu_{\bf s}-\mu_{\bf t}))^{i_k}
&\leq& C \|{\bf s}-{\bf t}\|_\infty^{\sum_{l=2}^n (l+1-\beta/d)i_l}\\
\label{eq:maj1}
&\leq& C\|{\bf s}-{\bf t}\|_\infty^{k+(1-\beta/d)(k/2)}
\end{eqnarray*}
where we use $\sum_{l=2}^kli_l=k$ and $\sum_{l=2}^ki_l\leq k/2$ together with $\|{\bf s}-{\bf t}\|_\infty\leq 1$. 
We deduce
$$ 
\ee\left[\left(\widetilde M_\rho(\mu_{\bf s}-\mu_{\bf t})\right)^k\right]
\leq C\|{\bf s}-{\bf t}\|_\infty^{(3-\beta/d)k/2}
$$
and using Corollary \ref{corol:moment},  
\begin{equation*}
\label{eq:maj2}
\ee\left[\left(J_a(\mu_{\bf s}-\mu_{\bf t})\right)^k\right]
\leq C\|{\bf s}-{\bf t}\|_\infty^{(3-\beta/d)k/2}.
\end{equation*} 
From the standard regularity criterion for random fields (see e.g. \cite[Th. 1.4.1]{Kun}), $J_a(\mu_{\bf t})$ has, almost surely, H\"older-continuous path for all index striclty less than $(3-\beta/d)/2-p/k$.
\CQFD


\subsection{Convergence in distributions space}
\label{sec:distribution}

The results on convergence in distributions space are based on Lemma \ref{lem:repres} below relating the generalized random field $\left(\widetilde M_\rho(\mu)\right)_{\mu\in\cD}$ and the parametric random field $\left(\widehat M_\rho({\bf t})\right)_{{\bf t}\in\real^d}$ defined by 
\begin{eqnarray*}
\widehat M_\rho(t)&=& \int_0^{t_1}\cdots \int_{0}^{t_d} \widetilde M_\rho(\delta_{\bf y})d{\bf y}\\
&=& \widetilde M_\rho(\mu_{\bf t}),\quad {\bf t}\in\real^d
\end{eqnarray*}
where $\mu_{\bf t}(d{\bf y})={\rm sign}(t_1)\cdots {\rm sign}(t_d)\ind_{[{\bf 0},{\bf t}]}({\bf y})d{\bf y}$ is the parametric family of signed uniform measures given in Example \ref{exa:1}. 
Recall that this family satisfies the property $(P_\gamma)$ for all $\gamma\geq 1$ so that Theorem \ref{theo:CV_functional1} indeed holds. 
We also define the continuous random fields on $\real^d$ which are the possible limits in ${\mathcal C}(\real^d)$ under, respectively, the large ball regime and the intermediate ball regime: 
$$
\widehat Z_\alpha({\bf t})=  Z_\alpha(\mu_{\bf t}), \quad \widehat J_a({\bf t})= J_a(\mu_{\bf t}),\quad {\bf t}\in\real^d.\\
$$
\begin{lemme}
\label{lem:repres}
For all $\phi\in \cD(\real^d)$, 
$$
\widetilde M_\rho(\phi({\bf y})d{\bf y})=(-1)^d \int_{\real^d} \widehat M_\rho({\bf y})\frac{\partial^d \phi}{\partial y_1\cdots\partial y_d}({\bf y})d{\bf y}.
$$
\end{lemme}
\begin{Proof}
This follows from successive integration by parts:
\begin{eqnarray*}
&&\widetilde M_\rho(\phi({\bf y})d{\bf y})\\
&=&\int_{\real^d} \widetilde M_\rho(\delta_{\bf y})\phi({\bf y})d{\bf y} \\
&=&(-1)^d \int_{\real^d} \left(\int_0^{y_1}\cdots \int_{0}^{y_d} \widetilde M_\rho(u_1, \dots, u_d)\ du_1\cdots du_d \right) \left(\frac{\partial^d \phi}{\partial y_1\cdots\partial y_d}({\bf y}) \right)d{\bf y}\\
&=&(-1)^d  \int_{\real^d} \widehat M_\rho({\bf y})\frac{\partial^d \phi}{\partial y_1\cdots\partial y_d}({\bf y})d{\bf y}. 
\end{eqnarray*}
\CQFD
\end{Proof}

We prove now Theorem \ref{theo:distrib}. 
The proof relies on Lemma \ref{lem:repres} and on Theorem \ref{theo:CV_functional1} stating 
that $\widehat M_\rho$ converge in ${\mathcal C}(\real^d)$ to $\widehat Z_\alpha$ (resp. $\widehat J_a$) in the large ball (resp. intermediate ball) regime.

\medskip
\noindent
{\bf Proof of Theorem \ref{theo:distrib}.}\\
{\bf First point:} The random field $\widetilde M_\rho$ is a bounded linear operator on $\cD(\real^d)$. 
Indeed, for $\phi\in\cD(\real^d)$ whose support is included in $[-T,T]^d$, Lemma \ref{lem:repres} entails:
$$ 
\left|\widetilde M_\rho(\phi({\bf y})d{\bf y})\right|\leq (2T)^d \left(\sup_{[-T,T]^d}|\widehat M_\rho| \right) \left(\sup_{[-T,T]^d} \frac{\partial^d \phi}{\partial y_1\cdots\partial y_d}\right).
$$
This proves the continuity of the linear application on $\cD(\real^d)$ and hence that $\widetilde M_\rho$ can be seen as a random distribution.\\
{\bf Second and third point:} 
Let $I:{\mathcal C}(\real^d)\to \cD'(\real^d)$ be the canonical injection given by 
$$
I(f):\phi\mapsto \int_{\real^d}f({\bf y})\phi({\bf y})d{\bf y}, \quad f\in {\mathcal C}(\real^d),
$$ 
and define the differential operator $D: \cD'(\real^d)\to \cD'(\real^d)$ by
$$
D(s):\phi \mapsto (-1)^d s\left(\frac{\partial^d \phi}{\partial y_1\cdots\partial y_d} \right),\quad s\in \cD'(\real^d).
$$
With these notations, Lemma \ref{lem:repres} states that $\widetilde M_\rho = (D\circ I)(\widehat M_\rho)$.
Since, the operators $D$ and $I$ are continuous, and $\widehat M_\rho$ converges in ${\mathcal C}(\real^d)$ as $\rho\to 0$ to $\widehat Z_\alpha$ (resp. $\widehat J_a$) in the large ball regime (resp. in the intermediate ball regime), 
the continuous mapping theorem implies that $\widetilde M_\rho$ weakly converges in $\cD'(\real^d)$ to $(D\circ I)(\widehat Z_\alpha)$ (resp. to $(D\circ I)(\widehat J_\alpha)$). 
Finally, since weak convergence in $\cD'(\real^d)$ implies {\it fdd} convergence on $\cD(\real^d)$, we have $(D\circ I)(\widehat Z_\alpha)\stackrel{{\it fdd}}{=}Z_\alpha$ and $(D\circ I)(\widehat J_a)\stackrel{{\it fdd}}{=}J_a$. This shows that $Z_\alpha$ and $J_a$ have modifications that are continuous on $\cD(\real^d)$.
\CQFD


\appendix

\section{Technical results}
\label{sec:technical}

\begin{lemme}
\label{lemme:muh}
Let $\gamma\geq 1$ and $\mu(d{\bf y})=\phi({\bf y})d{\bf y}$ with $\phi\in L^1(\real^d)\cap L^\gamma(\real^d)$
and $h\in L^1(\real^d)\cap L^\gamma(\real^d)$. 
Then
\begin{equation}
\label{eq:bound0}
\int_{\real^d} |\mu[\tau_{{\bf x},r}h]|^\gamma d{\bf x}\leq (r^d\|\phi\|_1^{\gamma}  \|h\|_\gamma^\gamma)\wedge (r^{\gamma d}\|\phi\|_\gamma^\gamma \|h\|_1^\gamma).
\end{equation}
\end{lemme}
\begin{Proof}
Since $h\in L^\gamma(\real^d)$ and $\mu(d{\bf y})=\phi({\bf y})d{\bf y}$ with $\phi\in L^1(\real^d)$, the H\"older inequality entails:   
\begin{eqnarray}
\nonumber
 \int_{\real^d} |\mu[\tau_{{\bf x},r}h ]|^\gamma d{\bf x}
&=& \int_{\real^d} \left|\int_{\real^d}h\left(\frac{{\bf y}-{\bf x}}{r}\right)\phi({\bf y})d{\bf y}\right|^\gamma d{\bf x}\\
\nonumber
&\leq& \|\phi\|_1^{\gamma-1}  \int_{\real^d\times\real^d} \left|h\left(\frac{{\bf y}-{\bf x}}{r}\right)\right|^\gamma |\phi({\bf y})|d{\bf y}  d{\bf x}\\
\nonumber 
&=& r^d\|\phi\|_1^{\gamma-1} \int_{\real^d\times\real^d} |h({\bf y})|^\gamma |\phi(ry+x)|d{\bf y}d{\bf x} \\
\label{eq:bound3}
&=& r^d\|\phi\|_1^\gamma\|h\|_\gamma^\gamma. 
\end{eqnarray}
On the other hand, still using H\"older inequality but with $\phi\in L^\gamma(\real^d)$ and $h\in L^1(\real^d)$, 
we have: 
\begin{eqnarray}
\nonumber 
\int_{\real^d} |\mu[\tau_{{\bf x},r}h]|^\gamma d{\bf x} 
&=& r^{\gamma d} \int_{\real^d} \left|\int_{\real^d}h({\bf y})\phi(ry+x)d{\bf y}\right|^\gamma d{\bf x}\\
\nonumber
&\leq &  r^{\gamma d}\int_{\real^d} \left(\int_{\real^d}|h({\bf y})|d{\bf y}\right)^{\gamma-1}\int_{\real^d}|\phi(ry+x)|^\gamma|h({\bf y})|d{\bf y} d{\bf x}\\
\label{eq:bound4}
&\leq &  r^{\gamma d}\|\phi\|_\gamma^\gamma \|h\|_1^\gamma. 
\end{eqnarray}
The bounds \eqref{eq:bound3} and \eqref{eq:bound4} together entail \eqref{eq:bound0}. 
\CQFD
\end{Proof}

\bigskip
\noindent
The following result proves the continuity required to apply Lemmas 2 and 3 instead of Lemma 6 in \cite{KLNS} in the modification in Section \ref{sec:fdd}. 

\begin{lemme}
\label{lem:cont}
Suppose that the fading function $h$ satisfies \eqref{eq:halpha}. 
For $\mu(d{\bf y})\in L^1(\real^d)\cap L^\alpha(\real^d)$, 
the application 
$r\mapsto \int_{\real^d}\Psi_G\left(\mu[\tau_{{\bf x},r}h]\right)d{\bf x}$ is continuous on $(0,+\infty)$. 
The same holds true for $\int_{\real^d}\Psi_\alpha\left(\mu[\tau_{{\bf x},r}h]\right)d{\bf x}$
with $\Psi_\alpha(\theta)=-\sigma^\alpha|\theta|^\alpha(1+i b \varepsilon(\theta)\tan(\pi\alpha/2))$
and $\sigma$, $b$ given in \eqref{eq:local}.
\end{lemme} 
\begin{Proof}
Observe first that, for  $r_n\to r_0>0$, $n\to\infty$, we have: 
\begin{equation}
\label{eq:continu}
\lim_{n\to\infty} \mu[\tau_{{\bf x},r_n}h]=\mu[\tau_{{\bf x},r_0}h]\quad d{\bf x}\mbox{-a.e.}
\end{equation} 
This is a standard application of Lebesgue's convergence theorem. 
Indeed, since $h$ is almost-everywhere continous, 
we have $\tau_{{\bf x},r_n}h({\bf y})\to \tau_{{\bf x},r_0}h({\bf y})$, $n\to+\infty$. 
From the definition of $h^\ast$, the convergence is bounded for all $n\geq 1$:
$$
|\tau_{{\bf x},r_n}h({\bf y})\phi({\bf y})| \leq \tau_{{\bf x},R} h^\ast(y)|\phi({\bf y})|,
\quad x,y\in\real^d
$$
with $R=\sup\{r_n;n\geq 1\}$ and the bound is $d{\bf x}$-integrable 
since Lemma \ref{lemme:muh} applied to $h^\ast\in L^1(\real^d)\cap L^\alpha(\real^d)$ rewrites
$$
\int_{\real^d} \left|\int_{\real^d}\tau_{{\bf x},R} h^\ast({\bf y})|\phi({\bf y})|d{\bf y}\right|^\alpha d{\bf x}<+\infty.
$$
A second application of Lebesgue's convergence theorem yields 
$$
\lim_{n\to+\infty} \int_{\real^d}\Psi_G(\mu(\tau_{{\bf x},r_n}h))d{\bf x}=\int_{\real^d}\Psi_G(\mu(\tau_{{\bf x},r_n}h))d{\bf x}.
$$
The convergence \eqref{eq:continu} together with the continuity of $\Psi_G$ imply indeed the pointwise convergence
$$
\lim_{n\to+\infty} \Psi_G(\mu(\tau_{{\bf x},r_n}h))=\Psi_G(\mu(\tau_{{\bf x},r_0}h))\quad d{\bf x}\mbox{-a.e.}
$$
The convergence is bounded by $C|\mu(\tau_{{\bf x},R}h^\ast)|^\alpha$ since $\Psi_G(u)\leq C|u|^\alpha$, 
which is integrable by Lemma \ref{lemme:muh} applied to $h^\ast\in L^1(\real^d)\cap L^\alpha(\real^d)$.
A similar proof holds for
$$
\lim_{n\to\infty} \int_{\real^d}\Psi_\alpha(\mu(\tau_{{\bf x},r_n}h))d{\bf x}=\int_{\real^d}\Psi_\alpha(\mu(\tau_{{\bf x},r_n}h))d{\bf x}.
$$
\CQFD
\end{Proof}

\medskip
\noindent
The following result estimates the truncated moments of $F$. 
In particular, the condition $\beta>d$ ensures that $F$ has a finite moment of order $d$.
\begin{lemme}
\label{lemme:F_truncated}
For $\delta>0$, when $u\to+\infty$, we have:
\begin{equation}
\label{eq:moment1}
\int_0^u r^\delta F(dr) \sim 
\left\{\begin{array}{ll}
\mbox{Cst}&\mbox{ if } \delta<\beta\\
\beta C_\beta \ln u& \mbox{ if } \delta=\beta\\
\frac{\beta}{\delta-\beta} C_\beta u^{\delta-\beta}& \mbox{ if } \delta>\beta
\end{array}
\right .
\end{equation}
and for $0<\delta<\beta$, when $u\to+\infty$, we have:
\begin{equation}
\label{eq:moment2}
\int_u^{+\infty} r^\delta F(dr) \sim \frac{\beta}{\beta-\delta}C_\beta u^{\delta-\beta}.
\end{equation}
Moreover when $\delta>\beta$, we have the global bound 
\begin{equation}
\label{eq:global}
\int_0^u r^\delta F(dr)\leq Cu^{\delta-\beta}.
\end{equation}
\end{lemme}
\begin{Proof}
Let $R$ be a random variable with distribution $F$. 
We have 
\begin{eqnarray}
\nonumber
\int_0^u r^\delta F(dr)&=&\ee\left[R^\delta\ind_{R\leq u}\right]
=\int_\Omega\int_{\real^+}\ind_{\{R\leq u\}}\ind_{\{t\leq R^\delta\}} dtdP\\ 
\nonumber
&=&\int_{\real^+} \pit(t^{1/\delta}\leq R\leq u) dt
=\int_0^{u^\delta}\pit(t^{1/\delta}\leq R\leq u) dt\\ 
\nonumber
&=&\delta\int_0^u\pit(s\leq R\leq u) s^{\delta-1}ds\\
\label{eq:tech1}
&=&\delta\int_0^u\pit(R\geq s) s^{\delta-1}ds-\delta\int_0^u\pit(R>u) s^{\delta-1}ds.
\end{eqnarray} 
But using condition \eqref{eq:F2}, we have:
\begin{eqnarray*}
\int_0^u\pit(R\geq s) s^{\delta-1}ds\sim 
\left\{\begin{array}{ll}
\mbox{Cst}&\mbox{ if } \delta<\beta\\
C_\beta \ln u& \mbox{ if } \delta=\beta\\
\frac{C_\beta}{\delta-\beta} u^{\delta-\beta}& \mbox{ if } \delta>\beta 
\end{array}
\right .
\end{eqnarray*}
and $\delta\int_0^u\pit(R>u) s^{\delta-1}ds\sim C_\beta u^{\delta-\beta}$
from which \eqref{eq:moment1} easily derives. 
Next, 
\begin{eqnarray*}
\int_u^{+\infty} r^\delta F(dr)&=&\ee[R^\delta\ind_{R\geq u}]
=\int_\Omega\int_{\real^+} \ind_{\{R\geq u\}}\ind_{\{t\leq R^\delta\}} dtdP\\ 
&=&\int_{\real^+} \pit(R\geq \max(t^{1/\delta},u)) dt
=\delta \int_{\real^+} \pit(R\geq \max(s,u)) s^{\delta-1} ds\\
&=&\delta \int_0^u \pit(R\geq u) s^{\delta-1} ds
+\delta \int_u^{+\infty} \pit(R\geq s) s^{\delta-1} ds\\
&\sim& C_\beta u^{\delta-\beta}+ \frac{\delta}{\beta-\delta}C_\beta u^{\delta-\beta}
=\frac{\beta}{\beta-\delta}C_\beta u^{\delta-\beta}
\end{eqnarray*} 
which is \eqref{eq:moment2}.
Finally, since $\pit(R\geq s)\leq 1$, \eqref{eq:tech1} entails $\int_0^u r^\delta F(dr)=O(u^\delta)$
so that together with \eqref{eq:moment1}, it is easy to derive \eqref{eq:global}. 
\CQFD
\end{Proof}

\begin{lemme}
\label{lemme:intmuh}
Let $\gamma>\beta/d$ and $\mu(d{\bf y})=\phi({\bf y})d{\bf y}$ with $\phi\in L^1(\real^d)\cap L^\gamma(\real^d)$ 
and $h\in L^1(\real^d)\cap L^\gamma(\real^d)$. Then, for any $\rho>0$,
\begin{equation*}
\label{eq:A3}
\int_{\real^d\times\real^+} |\mu[\tau_{{\bf x},r}h]|^\gamma d{\bf x}F_\rho(dr)\leq M \rho^\beta \frac{(\gamma-1)\beta d}{(\gamma d-\beta)(\beta-d)}
(\|\phi\|_\gamma \|h\|_1)^{\frac{\gamma(\beta-d)}{(\gamma-1)d}}(\|\phi\|_1\|h\|_\gamma)^{\frac{(\gamma d-\beta)\gamma}{(\gamma-1)d}}, 
\end{equation*}
where $M$ is a constant depending only on $F$.
\end{lemme}
\begin{Proof}
Using Lemma \ref{lemme:muh}, we have: 
\begin{eqnarray*}
&&\int_{\real^d\times\real^+} |\mu[\tau_{{\bf x},r}h]|^\gamma d{\bf x}F_\rho(dr)\\
&\leq&\int_{\real^+}(r^d\|\phi\|_1^\gamma  \|h\|_\gamma^\gamma)\wedge (r^{\gamma d}\|\phi\|_\gamma^\gamma \|h\|_1^\gamma)\ F_\rho(dr)\\
&=&\int_{\real^+}(\rho^dr^d\|\phi\|_1^{\gamma}  \|h\|_\gamma^\gamma)\wedge (\rho^{\gamma d}r^{\gamma d}\|\phi\|_\gamma^\gamma \|h\|_1^\gamma)\ F(dr)\\
&=&\|\phi\|_\gamma^\gamma \|h\|_1^\gamma\rho^{\gamma d}\int_0^{c/\rho} r^{\gamma d}F(dr)+ \|\phi\|_1^{\gamma}  \|h\|_\gamma^\gamma\rho^{d}\int_{c/\rho}^{+\infty}r^d \ F(dr)
\end{eqnarray*}
with $c=((\|\phi\|_1 \|h\|_\gamma)/(\|\phi\|_\gamma \|h\|_1))^{\frac{\gamma}{(\gamma-1)d}}$.
Finally using the bound on the truncated moments of $F$ in Lemma \ref{lemme:F_truncated} 
with $\delta:=\gamma d>\beta$ in \eqref{eq:moment1} and $\delta:=d$ in \eqref{eq:moment2}, 
we derive 
\begin{equation}
\nonumber
\int_{\real^d\times\real^+} |\mu[\tau_{{\bf x},r}h]|^\gamma d{\bf x}F_\rho(dr)\leq M\rho^\beta \left[\frac{\beta}{\gamma d-\beta}\|\phi\|_\gamma^\gamma \|h\|_1^\gamma c^{\gamma d-\beta}+\frac{\beta}{\beta-d}\|\phi\|_1^{\gamma}  \|h\|_\gamma^\gamma c^{d-\beta}\right]
\end{equation}
where $M\in(0,+\infty)$ depends only on $F$. 
The result is obtained, after cancellation, by replacing $c$ by its definition. 
\CQFD
\end{Proof}

\medskip
\noindent
The next result collects explicit formulas for the characteristic function 
and the cumulants of the rescaled and centered random variable $\widetilde M_\rho(\mu)$.
It is based on standard results for Poisson integrals, see \cite{Kallenberg}. 
\begin{lemme}
\label{prop:cumulant}
\begin{enumerate}
\item The characteristic function of  $\widetilde M_\rho(\mu)$ writes:
\begin{equation*}
\varphi_{\widetilde M_\rho(\mu)}(\theta)
=\exp\left(\int_{\real^d\times\real^+\times\real} \Psi_G\left(n(\rho)^{-1}\theta \mu[\tau_{{\bf x},r}h ]\right)\lambda(\rho)d{\bf x}F_\rho(dr)\right)
\end{equation*}
where $\Psi_G(u)=\int_\real (e^{ium}-1-ium )G(dm)$.
\item Suppose $G$ has a finite moment of order $k\geq 1$, $\mu(d{\bf y})=\phi({\bf y})d{\bf y}$ with $\phi\in L^1(\real^d)\cap L^k(\real^d)$
and $h\in L^1(\real^d)\cap L^k(\real^d)$. 
Then $\widetilde M_\rho(\mu)$ has a finite moment of order $k$ and its $k$ first cumulants are given by:
\begin{eqnarray*}
c_1(\widetilde M_\rho(\mu))\!\!\!\!&=&\!\!\!\!0,\\
c_l(\widetilde M_\rho(\mu))\!\!\!\!&=&\!\!\!\!\frac{\lambda(\rho)}{{n(\rho)}^l}\left(\int_{\real}m^l\ G(dm)\right) \left(\int_{\real^d\times\real^+} (\mu[\tau_{{\bf x},r}h])^l d{\bf x}\ F_\rho(dr)\right),\ 2\leq l\leq k.
\end{eqnarray*}
\end{enumerate}
\end{lemme}
Note that the finiteness of 
$\int_{\real^d\times\real^+} |\mu[\tau_{{\bf x},r}h]|^l\ d{\bf x}F_\rho(dr)$
comes from Lemma \ref{lemme:muh} when $\mu\in L^1(\real^d)\cap L^k(\real^d)$ and $h\in L^1(\real^d)\cap L^k(\real^d)$. 

\bibliographystyle{alea2}
\bibliography{}

\end{document}